\documentclass[11pt,reqno]{amsart}
\usepackage{hyperref}
\usepackage{enumerate,fullpage}
\usepackage{amssymb,amsmath,graphicx,amsthm}
\usepackage{color}
\usepackage{soul}

\newtheorem{theorem}{Theorem}

\newtheorem{lemma}[theorem]{Lemma}

\theoremstyle{remark}
\newtheorem{remark}{Remark}
\theoremstyle{definition}

\newcommand{\la}{\langle}
\newcommand{\ra}{\rangle}

\newcommand{\no}[1]{\|{#1}\|}
\newcommand{\Om}{\Omega}
\newcommand{\om}{\omega}

\newcommand{\we}{\wedge}
\newcommand{\p}{\partial}

\newcommand{\z}{\bar z}
\newcommand{\w}{\bar w}

\def\di{\partial}
\def\dib{\bar\partial}
\def\Im{\T{Im}\,}
\def\Re{\T{Re}\,}
\def\T{\text}
\def\C{\mathbb C}
\def\R{\mathbb R}
\def\N{\mathbb N}

\def\p{\partial}
\numberwithin{equation}{section}
\renewcommand\Re{\operatorname{Re}}
\renewcommand\Im{\operatorname{Im}}

\def\T{\text}

\DeclareMathOperator{\supp}{supp}

\DeclareMathOperator{\dist}{dist}

\DeclareMathOperator{\Rre}{Re}

\newcommand{\dbar}{\bar\partial}

\newcommand{\dbarb}{\bar\partial_b}

\newcommand{\vp}{\varphi}

\newcommand{\bd}{b}

\newcommand{\loc}{\text{loc}}

\newcommand{\ep}{\epsilon}

\begin{document}

	\title[Zero sets of functions of the Nevanlinna class]{Zero sets of functions in the Nevanlinna class and the $\bar\partial_b$-equation on convex domains of general type in $\C^2$}
	

\author[T.~V.~Khanh, and A.~Raich]{Tran Vu Khanh and Andrew Raich}
\address{T.V.~Khanh}
\address{School of Mathematics and Applied Statistics, University of Wollongong, NSW, Australia,  2522}
	\email{tkhanh@uow.edu.au}
		
\address{A. Raich}
\address{Department of Mathematical Sciences, SCEN 327, 1 University
  of Arkansas, Fayetteville, AR 72701} 
\email{araich@uark.edu}

\thanks{The first author was partially supported by Australian Research Council DE160100173.
The second author was partially supported by NSF grant DMS-1405100. He also gratefully acknowledges the support provided by
Wolfson College at the University of Cambridge during his Fellowship year there.}


\maketitle

\begin{abstract}
The purpose of this paper is to characterize the zero sets of holomorphic functions in the Nevanlinna class on a class of convex domains of infinite type in $\C^2$. Moreover, we also obtain $L^p$ estimates, 
$1 \leq p \leq \infty$,  for a particular solution of the tangential Cauchy-Riemann equation on the boundaries of these domains. 
\end{abstract}

\section{Introduction}\label{sec:intro}
\subsection{Background}
Let $\Om$ be a $C^\infty$-smooth, bounded domain with smooth defining function $\rho$. In several complex variables, characterizing the zero sets of holomorphic functions in the Nevanlinna class is
closely related to the Poincar\'e-Lelong equation
\[
i \p\dbar u = \alpha \quad\text{on }\Om
\]
for $d$-closed, smooth $(1,1)$-forms $\alpha$. Not surprisingly, solving the Poincar\'e-Lelong equation with estimates often amounts to studying the Cauchy-Riemann equation
\[
\dbar u =f \quad\text{on }\Om
\]
where $f$ is a $\dbar$-closed $(0,1)$-form on $f$. For investigating the Nevanlinna class, the type of estimates that are useful are solutions $u$ with boundary values in $L^1(\bd\Om)$. The techinque that we use to solve the
boundary value estimate in turn yields a solution to the tangential Cauchy-Riemann equation
\[
\dbarb u =f
\]
in $L^p(\bd\Om)$, $1 \leq p \leq \infty$ where $f$ is $\dbarb$-closed $(0,1)$. In fact, our technique is produces a gain in $L^\infty$ and maps $L^\infty$ to an appropriate $f$-H\"older space. Therefore, our results
concern the related problems of complex varieties that are zero sets of Nevanlinna functions, the Cauchy-Riemann, and tangential Cauchy-Riemann equations.

It is well-known that if $h\in N(\Om)$ then the zero divisor $X_h$ of $h$ satisfies the Blaschke condition. Whether or not the converse is true, namely,
{\it ``If $h\in H(\Om)$ and $X_h$ satisfies the Blaschke condition, does there exists $f \in N(\Om)$ so that $X_f = X_h$?" } has been extensively studied over the past forty years.
The $n=1$ case is a classic one-variable result in the complex plane. In contrast, for $n\geq 2$, the Blaschke condition for a divisor no longer suffices to be the zero set of a Nevanlinna function or even
the zero set of an $H^p$ function \cite{Var80}. 
There are cases, however, where the Blaschke condition is sufficient. Namely, the sufficiency is known when $\Om$ is
\begin{itemize}
  \item a strongly pseudoconvex domain \cite{Gru75, Sko76};
  \item a pseudoconvex domain in $\C^2$  of finite type \cite{ChNaSt92,Sha89}; 
  \item a complex or real ellipsoid by \cite{BoCh82, Sha91h};   
  \item a convex domain of strictly finite type and/or finite type \cite{BrChDu98, ChDuMo14, Cum01n,DiMa01}. 
\end{itemize}
Furthermore, the positive answer still holds on the following infinite type example:
\begin{equation}\label{Dalpha}
D_\alpha=\Big\{(z_1,z_2)\in \C^2: |z_1|^2+\exp\left(1+\frac{2}{\alpha}-\frac{1}{|z_2|^\alpha}\right)<1\Big\}
\end{equation}
 with $0<\alpha<1$ (see \cite{AhCh02}).

In this paper, we shall prove the converse is true for the class of convex domains in $\C^2$ of general type by Khanh in \cite{Kha13, HaKhRa14}. 
All known examples of convex domains in $\C^2$ are covered by our class.

Establishing $L^p$ and H\"older estimates for solutions of the $\dib_b$-equation is a fundamental question in several complex variables. 
It has been extensively investigated on classes of domains of finite type such as strongly pseudoconvex domains \cite{FoSt74e,Hen77a,Hen77b}, convex domains \cite{Ale05, Sha91h}, domains with a diagonalizable Levi form
\cite{FeKoMa90},
domains where the Levi form has comparable eigenvalues \cite{Koe02},
decoupled domains \cite{NaSt06}, and pseudoconvex domains in $\C^2$ \cite{Chr88, FeKo88}. See also \cite{Wu98, LTSh05}.
 Theorem~\ref{t3} provides the first example of $L^p$, $1 \leq p \leq \infty$, and H\"older estimates on an infinite type domain.

\subsection{The class of general type convex domains}
Our setup is the following: $\Omega\subset\C^2$ is a smooth, bounded domain. For 
each $p\in\bd\Om$, the curvature of $\bd\Om$ at $p$ is captured by local coordinates $z_p = T_p(z)$ where $T_p$ is a $\C$-linear 
transformation that sends $p$ to the origin. Additionally, there exist a global
defining function $\rho$ and functions $F_p$ and $r_p$ satisfying
\begin{equation}\label{OmC}
\Om_p=T_p(\Om)=\{z_p=(z_{p,1},z_{p,2})\in \C^2:\rho(T_p^{-1}(z_p))=F_p(|z_{p,1}|^2)+r_p(z_p)<0\}
\end{equation}
or 
\begin{equation}\label{OmR}
\Om_p=T_p(\Om)=\{z_p=(z_{p,1},z_{p,2})\in \C^2: \rho(T^{-1}_p(z_p))=F_p(x_{p,1}^2)+r_p(z_p)<0\}
\end{equation}
where $z_{p,j}=x_{p,j}+iy_{p,j}$, $x_{p,j}, y_{p,j}\in\mathbb{R}$, $j=1,2$, and $i=\sqrt{-1}.$ We also require that the functions $F_p:\R\to\R$ and $r_p:\C^2\to\R$ satisfy:
\begin{enumerate}[i.] 
\item $F_p(0)=0$;
\item $F_p'(t), F_p''(t), F_p'''(t)$, and $\left(\dfrac{F_p(t)}{t}\right)'$ are nonnegative on $[0,\tilde d^2_p)$ for some $\tilde d_p>d_p$;
\item $r_p(0)=0$ and $\dfrac{\di r_p}{\di z_{p,2}}\ne 0$ on $\bd\Om$ with $|z_{1,p}|\le \delta$;
\item $r_p$ is convex,
\end{enumerate}
where $d_p$ is the diameter of $\Om_p$ and $\delta$ is a small number independent of $p$.\\

This class of domains includes several well-known examples. 
If $\Omega$ is of finite type $2m$, then $F_p(t) = t^m$ at the points of type $2m$.
On the other hand, if $F_p(t)=\exp(-1/t^{\alpha})$, then $\Omega$ is of infinite type at $p$, and this is our main case of interest. 
 Our hypotheses include the following three classes of infinite type domains: the complex ellipsoid
\begin{eqnarray}\label{Om1}
\Om=\bigg\{z=(z_1,z_2)\in \C^2:\sum_{j=1}^2 \exp\left(-\frac{1}{|z_j|^{\alpha_j}} \right)\le e^{-1}\bigg\};
\end{eqnarray} 
the real ellipsoid
\begin{eqnarray}\label{Om2}
 \Om=\bigg\{z = (x_1+iy_1,x_2+iy_2)\in\C^2:\sum_{j=1}^2\exp\left(-\frac{1}{|x_j|^{\alpha_j}}\right)+\exp\left(-\frac{1}{|y_j|^{\beta_j}}\right)\le e^{-1}\bigg\};
\end{eqnarray}
and  the mixed case
\begin{eqnarray}\label{Om3}
\Om=\bigg\{z = (x_1+iy_1,x_2+iy_2)\in\C^2:\exp\left(-\frac{1}{|x_1|^{\alpha_1}}\right)+\exp\left(-\frac{1}{|y_1|^{\beta_1}}\right)+\exp\left(-\frac{1}{|z_2|^{\alpha_2}} \right)\le e^{-1}\bigg\}
\end{eqnarray}
where $\alpha_j, \beta_j>0$. Moreover, our setting also includes  a tube domain of infinite type at $0$
\begin{eqnarray}\label{Om4}
\Om=\{z = (x_1+iy_1,x_2+iy_2)\in\C^2:\exp\left(1-\frac{1}{|x_1|^{\alpha_1}}\right)+\chi(y_1)+|z_2|^2\le 1\}
\end{eqnarray}
where $\chi$ is a convex function and $\chi(y_1)=0$ when $|y_1|<\delta$ and $\alpha_1>0$ for $j=1,2$.

\subsection{Notation}
For an excellent discussion of the Nevanlinna class, complex varieties, (positive) currents, and (irreducible) divisors, we strongly encourage the reader to consult Range \cite{Ran86} and Noguchi and Ochiai \cite{NoOc90}.

The \emph{Nevanlinna class} for $\Om$, denoted by $N(\Om)$, is defined by 
$$N(\Om)=\left\{h\in H(\Om):\sup_{\epsilon>0} \int_{\bd\Om_\epsilon}\left|\log|h(z)|\right|d\sigma_{\bd\Om_\epsilon}(z)<+\infty\right\},$$
where $H(\Om)$ is the space of holomorphic functions on $\Om$, $\bd\Om_\epsilon=\{z:\rho(z)=-\epsilon\}$ for small $\epsilon>0$, 
and $d\sigma_{\bd\Om_\epsilon}(z)$ denotes the Euclidean surface measure on $\bd\Om\epsilon$. If $X\subset \Om$ is a complex variety with irreducible decomposition 
$$X=\bigcup_kX_k$$
and $n_k\in\N$ are positive integers for each $k$, the divisor $\hat X:=\{X_k,n_k\}$ is said to satisfy the \emph{Blaschke condition} if 
$$\sum_k n_k\int_{X_k}|\rho(z)|d\mu_{X_k}(z)<\infty$$
where $d\mu_{X_k}$ is the induced surface area measure on $X_k$.

\subsection{Main results}
We have three main results.
\begin{theorem}\label{main1} Let $\Om$ be a bounded domain in $\C^2$. Assume that for any $p\in \bd\Om$
\begin{enumerate}[i.] 
\item $\Om$ is defined by (\ref{OmC}) and $\int_0^{d_p}| \log F_p(t^2)|\,dt<\infty$ for all $p\in \bd\Omega$, or
\item $\Om$ is defined by (\ref{OmR}) and $\int_0^{d_p}|\log (t) \log F_p(t^2)|\,dt<\infty$ for all $p\in \bd\Omega$.
\end{enumerate}
Then for any divisors $\hat X$ in $\Om$ satisfying the Blaschke condition there is a function $h\in N(\Om)$ such that $\hat X$ is the zero divisor of $h$.
\end{theorem} 
To apply Theorem~\ref{main1} to the domains $\Om$ defined by any of \eqref{Om1} -  \eqref{Om4},  we are forced to require that $\alpha_j<1$, though any $\beta_j>0$ is permissible.  
The crucial step to prove  Theorem~\ref{main1} is the existence a solution of the $\dib$-equation that satisfies both of the following conditions:  the solution is
(i) smooth if data is smooth, and (ii) bounded in $L^1(\bd\Om)$. 
\begin{theorem}\label{t2} Let $\Omega$ be a domain satisfying the hypotheses of Theorem~\ref{main1}. 
Then for any $\dib$-closed, smooth $(0,1)$-form $\phi$  on $\Om$ so that $\no{\phi}_{L^1(\bd\Om)}<\infty$,  there exists a smooth function $u$ such that 
\begin{eqnarray}\label{dib}\dib u=\phi \quad\T{on $\Om$}\end{eqnarray}
and  
\begin{eqnarray}\label{L1bOm}
\no{u}_{L^1(\bd\Om)}\le c\no{\phi}_{L^1(\bd\Om)}
\end{eqnarray}
where $c>0$ is independent of $\phi$. 
\end{theorem}
\begin{remark}\label{rm1}
	The constant $c$ in \eqref{L1bOm} depends on the geometric type and diameter of $\Om$. It is, however,  uniformly bounded if $\Om$ satisfies the hypothesis of Theorem~\ref{main1} and diameter of $\Om$ is bounded. 
\end{remark}
Let $u$ be the Henkin solution to \eqref{dib} given by \eqref{Henkin} below. Then the smoothness of $u$, given the smoothness of $\phi$,  
is a consequence of Theorem  3 in \cite{Ran92}. Therefore we only need to prove that the inequality \eqref{L1bOm} holds for this $u$. 
The proof will be given in Section 2. A small modification of the technique to prove
\eqref{L1bOm} yields $L^p$-estimates for the tangential Cauchy-Riemann equation, a significant and new result in its own right.
\begin{theorem}\label{t3} Let $\Omega$ be a domain satisfying the hypotheses of Theorem~\ref{main1} and $p\in [1,\infty]$. Let $\phi$ be a $(0,1)$-form in $L^p(\bd\Om)$, satisfying the compatibility condition $\int_{\bd\Om}\phi\wedge \alpha=0$ for every continuous up the boundary, $\dib$-closed $(2,0)$-form $\alpha$ on $\Om$. Then  there exists  a function $u$ on $\bd\Om$ such that 
\begin{eqnarray}\label{dib-b}\dib_b u=\phi \quad\T{ on $\bd\Om$ }
\end{eqnarray}
and 
\begin{eqnarray}\label{LpbOm}
\no{u}_{L^p(\bd\Om)}\le c\no{\phi}_{L^p(\bd\Om)}
\end{eqnarray}
where $c>0$ is independent of $\phi$.\\

Moreover, in the case $p=\infty$, we obtain a ``gain" for the solution of $\dib_b$ into the $f$-H\"older spaces, that is, 
\begin{equation}\label{Linftyk}
\|u\|_{\Lambda^f(\bd\Om)} \leq c \|\phi\|_{L^\infty(\bd\Om)}.
\end{equation}
The function $f$ is defined by $f(d^{-1}):=\inf_{p\in \bd\Om}\left(\displaystyle\int_{0}^{d}\frac{\sqrt{F^*_p(t)}}{t}dt\right)^{-1}$
when $\Om$ is defined by \eqref{OmC} and by $f(d^{-1}):=\inf_{p\in \bd\Om} \left(\displaystyle\int_{0}^{d}\frac{\sqrt{F_p^*(t)}|\log \sqrt{F^*(t)}|}{t}dt\right)^{-1}$ when $\Om$ is defined by \eqref{OmR}. Here the superscript $^*$ denotes the inverse function and the $f$-H\"older space $\Lambda^f$ is defined by 
\begin{equation}
\Lambda^f(\bd\Omega)=\left\{u : ||u||_{\Lambda^f(\bd\Om)}:=||u||_{L^\infty(\bd\Om)}+\sup_{X(t)\in \mathfrak C, 0\le t\le1}\{f(|t|^{-1})|u(X(t))-u(X(0))|\}<\infty\right\}
\end{equation}
where $\mathfrak C=\{X(t):t\in [0,1]\to X(t)\in \bd\Om \T{ is $C^1$ and $|X'(t)|\le 1$}\}$.
\end{theorem}

\section{Proof of Theorem \ref{t2} and Theorem~\ref{t3}}\label{sec:proof of Thm 2,3}

We first assume that the origin is in $\bd\Om$, the functions $F=F_0$ and $r=r_0$ satisfy conditions (i)-(iv) from  Section \ref{sec:intro}, and
\begin{equation}\label{OmC0}
\Om=\{z=(z_{1},z_{2})\in \C^2:\rho(z)=F(|z_{1}|^2)+r(z)<0\}
\end{equation}
or 
\begin{equation}\label{OmR0}
\Om=\{z=(z_{1},z_{2})\in \C^2:\rho(z)=F(|x_{1}|^2)+r(z)<0\}
\end{equation}
where $z_{j}=x_{j}+iy_{j}$, $x_{j}, y_{j}\in\mathbb{R}$, $j=1,2$. Here we only need to consider $F'(t^2)\not=0$, for otherwise $\Om$ is strictly convex, and the proof of Theorem~\ref{t2} and \ref{t3} is known. Let
the support function for $\Om$ be defined by 
$$\Phi(\zeta,z)=\sum_{j=1}^2 \frac{\di \rho(\zeta)}{\di \zeta_j }(\zeta_j-z_j).$$
We are going to estimate $\Re\{\Phi(\zeta,z)\}$ for $\zeta,z$ in a neighborhood of $\bd\Om$.
 
 \noindent{\bf Setting 1: $\Om$ is defined by \eqref{OmC0}.} The convexity of $r$ yields that 
\begin{eqnarray}\label{RePhi}\begin{split}
2\Re\{\Phi(\zeta,z)\}\ge \rho(\zeta)-\rho(z)+ F'(|\zeta_1|^2)|z_1-\zeta_1|^2+\Big[F(|z_1|^2)-F(|\zeta_1|^2)-F'(|\zeta_1|^2)(|z_1|^2-|\zeta_1|^2)\Big]
\end{split}\end{eqnarray}
for any  $\zeta,z$ in a neighborhood of  $\bd\Om$. 
\begin{lemma}\label{lm22} Let $\Om$ be defined by \eqref{OmC0} and $F$ satisfy both conditions (i)-(iv) from Section \ref{sec:intro} and 
$F'(0)=0$.  Let $\zeta,z$ be in a neighborhood of $\bd\Om$ and satisfy $\rho(\zeta)-\rho(z)\ge 0$.
	\begin{itemize}
		\item If $|\zeta_1|\ge |z_1-\zeta_1|$ then 
		$$	|\Phi(\zeta,z)|^k|z-\zeta|\gtrsim  \big[|\Im\Phi(\zeta,z)|+\rho(\zeta)-\rho(z)+F(|z_1-\zeta_1|^2)\big]^k|z_1-\zeta_1|.$$
		\item Otherwise, if $|\zeta_1|\le |z_1-\zeta_1|$, then
		\begin{eqnarray}
		\begin{split}
|\Phi(\zeta,z)|^k|z-\zeta|\gtrsim&\big[|\Im\Phi(\zeta,z)|+\rho(\zeta)-\rho(z)+F(\frac{1}{2}|\zeta_1|^2)\big]^k|\zeta_1|\\
&+\big[|\Im\Phi(\zeta,z)|+\rho(\zeta)-\rho(z)+F(\frac{1}{2}|z_1|^2)\big]^k|z_1|
		\end{split}
		\end{eqnarray}	
	\end{itemize}
	for $k=1,2$.
\end{lemma}
\begin{proof}
 The term in $[\cdots]$ in \eqref{RePhi} is nonnegative for any $\zeta,z$ near $\bd\Om$, so the hypothesis on the sizes of  $|\zeta_1|$ and $|z_1-\zeta_1|$ allow us to obtain
 \begin{equation}\label{a1}
 |\Phi(\zeta,z)|^k|z-\zeta| \gtrsim \begin{cases}(|\Im\Phi(\zeta,z)|+\rho(\zeta)-\rho(z)+F(|z_1-\zeta_1|^2))^k|z_1-\zeta_1| &\T{if~~}~~|\zeta_1|\ge |z_1-\zeta_1|\vspace{.2in},\\
 (|\Im\Phi(\zeta,z)|+\rho(\zeta)-\rho(z)+F(|\zeta_1|^2))^k|\zeta_1|. &\T{if~~}~~|\zeta_1|\le |z_1-\zeta_1|.\\
 \end{cases}
 \end{equation}
The remaining estimate to show is
\begin{equation}\label{2.6k}
|\Phi(\zeta,z)|^k|z-\zeta|\gtrsim(|\Im\Phi(\zeta,z)|+\rho(\zeta)-\rho(z)+F(\frac{1}{2}|z_1|^2))^k|z_1| 
\end{equation}
in the case $|\zeta_1|\le |z_1-\zeta_1|$. It
can be obtained using the argument of Lemma 3.2 in \cite{HaKhRa14}. 
For the reader's convenience, we outline the proof here.  
Start by comparing the relative sizes of $|\zeta_1|$ and $\frac{1}{\sqrt{2}}|z_1|$.
If $|\zeta_1|\ge \frac{1}{\sqrt{2}}|z_1|$, then the argument follows from the second line of \eqref{a1}. 
Otherwise, $|\zeta_1|\le \frac{1}{\sqrt{2}}|z_1|$, and this inequality implies both $|z_1|\ge |\zeta_1|$ and $|z_1-\zeta_1|\ge (1-\frac{1}{\sqrt{2}})|z_1|$. 
We then estimate the $[\dots]$ in \eqref{RePhi} by
\begin{eqnarray}
\begin{split}
[\dots]:=F(|z_1|^2)-F(|\zeta_1|^2)-F'(|\zeta_1|^2)(|z_1|^2-|\zeta_1|^2)
\ge F(|z_1|^2-|\zeta_1|^2)\ge F(\frac{1}{2}|z_1|^2), 
\end{split}
\end{eqnarray}
where the inequality uses the facts that $F'(0)=0$ and $F''$ is nondecreasing (see \cite[Lemma 4]{FoLeZh11} or \cite[Lemma 3.1]{HaKhRa14} for details). This completes the proof.
\end{proof}
 \noindent{\bf Setting 2: $\Om$ is defined by \eqref{OmR0}.} An argument analogous to that of Lemma~\ref{lm22} produces the following lemma.
\begin{lemma}\label{lm23} Let $\Om$ be defined by \eqref{OmR0} and  $F$ satisfy both conditions (i)-(iv) from Section \ref{sec:intro} and 
$F'(0)=0$. Suppose $\zeta,z$ are in a neighborhood of $\bd\Om$  and $\rho(\zeta)-\rho(z)\ge 0$. 
	\begin{itemize}
		\item If $|\xi_1|\ge |x_1-\xi_1|$, then 
		$$	|\Phi(\zeta,z)|^k|z-\zeta|\gtrsim |\Im\Phi(\zeta,z)|+\rho(\zeta)-\rho(z)+F(|x_1-\xi_1|^2))^k(|x_1-\xi_1|+|y_1-\eta_1|).$$
		\item Otherwise, if $|\xi_1|\le |x_1-\xi_1|$, then 
		\begin{eqnarray}
		\begin{split}
		|\Phi(\zeta,z)|^k|z-\zeta|\gtrsim&(|\Im\Phi(\zeta,z)|+\rho(\zeta)-\rho(z)+F(\frac{1}{2}|\xi_1|^2))^k(|\xi_1|+|y_1-\eta_1|)\\
		&+(|\Im\Phi(\zeta,z)|+\rho(\zeta)-\rho(z)+F(\frac{1}{2}|x_1|^2))^k(|x_1|+|y_1-\eta_1|)
		\end{split}
		\end{eqnarray}	
	\end{itemize}
	for $k=1,2$, where $z_1=x_1+iy_1$ and $\zeta_1=\xi_1+i\eta_1$.
\end{lemma}

	Both the proofs of Theorem \ref{t2} and \ref{t3} use the supporting function estimates of Lemma~\ref{lm22} and \ref{lm23}. 
	We first give the full proof of Theorem~\ref{t3} since it is more dedicate and indicate the changes necessary to prove Theorem~\ref{t2}.

\subsection{Proof of Theorem~\ref{t3}}
 \begin{proof}[Proof of Theorem~\ref{t3}]

 	Let $\phi$ be a $(0,1)$-form satisfying the compatibility condition, that is, 
	$\int_{\bd\Om}\phi\we \alpha=0$ for every continuous up to the boundary $\dib$-closed $(2,0)$-form $\alpha$ on $\Om$. M.C. Shaw \cite{Sha91h, Sha89}  showed that 
 	\begin{equation}\label{Tb}
 	u:=T_b\phi=H^+\phi-H^-\phi
 	\end{equation}
 	is an integral solution (in the distribution sense) to the $\dib_b$-equation, $\dib_bu=\phi$, on $\bd\Om$ where    
 	$$H^+\phi(z):=\frac{1}{4\pi^2}\lim_{\epsilon\to 0^+}\int_{\bd\Om}H(\zeta,z-\epsilon \nu(z))\phi(\zeta)\we \om(\zeta),\quad\quad H^-\phi(z):=\frac{1}{4\pi^2}\lim_{\epsilon\to 0^+}\int_{\bd\Om}H(z+\epsilon \nu(z),\zeta)\phi(\zeta)\we \om(\zeta).$$
 	Here $\nu(z)$ is the outward unit normal vector at $z\in \bd\Om$; $\om(\zeta)=d\zeta_1 \wedge d\zeta_2$; and $H(\zeta,z)$ is given by 
 	\begin{equation} \label{H}
 	H(\zeta,z)=\frac{\frac{\p \rho(\zeta)}{\p\zeta_1}(\bar\zeta_2-\bar z_2)-\frac{\p \rho(\zeta)}{\p\zeta_2}(\bar\zeta_1-\bar z_1)}{\Phi(\zeta,z)|\zeta-z|^2}.
 	\end{equation}
In order to prove \eqref{LpbOm}, we will prove that $\no{u}_{L^p(\bd\Om)}\lesssim \no{\phi}_{L^p(\bd\Om)}$ only for $p=1$ and $p=\infty$; then using Riesz-Thorin Interpolation Theorem we obtain $L^p$ estimates for all $p\in [1,\infty]$. \\

{\it Part I: Proof of $\no{u}_{L^1(\bd\Om)}\lesssim \no{\phi}_{L^1(\bd\Om)}$.}
Denote by 
$$u_\epsilon(z):=\frac{1}{4\pi^2}\int_{\bd\Om}H(\zeta,z-\epsilon \nu(z))\phi(\zeta)\we \om(\zeta)-\frac{1}{4\pi^2}\int_{\bd\Om}H(z+\epsilon \nu(z),\zeta)\phi(\zeta)\we \om(\zeta), \quad z\in \bd\Om. $$
It follows that $\lim_{\epsilon\to 0^+} u_\epsilon=u$ a.e. 
For $\phi\in L^1(\bd\Om)$ and we need to prove that $\no{u_\epsilon}_{L^1(\bd\Om)}\le \no{\phi}_{L^1(\bd\Om)}$ uniformly for small $\epsilon>0$. It then follows that $u_\ep \to u$ in $L^1(\bd\Om)$ by the Dominated
Convergence Theorem. We observe 
\begin{align*}
\no{u_\epsilon}_{L^1(\bd\Om)}=&\frac{1}{4\pi^2}\int_{z\in \bd\Om}\left|\int_{\zeta\in \bd\Om}\left(H(\zeta,z-\epsilon \nu(z)-H(z+\epsilon \nu(z),\zeta)\right)\phi(\zeta)\we \om(\zeta)\right|dS(z)\\
\lesssim& \int_{z\in \bd\Om}\int_{\zeta \in \bd\Om}|H(\zeta,z-\epsilon \nu(z))||\phi(\zeta)|dS(\zeta)dS(z)+\int_{z\in \bd\Om}\int_{\zeta\in \bd\Om}|H(z-\epsilon \nu(z),\zeta)||\phi(\zeta)|dS(\zeta)dS(z)\\
\lesssim&  \int_{z\in \bd\Om_\epsilon}\int_{\zeta \in \bd\Om}|H(\zeta,z)||\phi(\zeta)|dS(\zeta)dS(z)+\int_{z\in \bd\Om^\epsilon}\int_{\zeta\in \bd\Om}|H(z,\zeta)||\phi(\zeta)|dS(\zeta)dS(z)
\end{align*}
where $\Om_\epsilon=\{z\in \Om:\rho(z)< -\epsilon\}$ and $\Om^\epsilon=\{z\in \C^2\setminus\bar\Om:\rho(z)> \epsilon\}$.
As a consequence of Tonelli's Theorem, it suffices to prove that
\begin{equation}\label{new1}\begin{split}\iint_{(\zeta,z)\in \bd\Om\times \bd\Om_\epsilon}\big| H(\zeta,z)||\phi(\zeta)\big|\, dS(\zeta,z)+\iint_{(z,\zeta)\in \bd\Om^\epsilon \times \bd\Om }\big| H(z,\zeta)||\phi(\zeta)\big|\, dS(z,\zeta )\lesssim \no{\phi}_{L^1(\bd\Om)}<\infty.
\end{split}\end{equation}

Since $\bd\Om$ is compact, for any  $\delta>0$, there exist points $p_1,\dots,p_N\in \bd\Om$ so that $\bd\Om$ is covered by$\{B(p_j, \delta)\}_{j=1}^N$. 
After changing coordinates with the linear transformation $T_{p_j}$ as in Section~1 (keeping in mind $T_{p_j}(p_j)=0$),  we may assume the goal is to prove
\begin{eqnarray}\label{new2}\begin{split}
&\underset{I}{\underbrace{\iint_{(\zeta_{p_j},z_{p_j})\in (B(0,\delta)\cap (\bd\Om_{p_j})_\ep)\times \bd\Om_{p_j}}\Big |H\big(T^{-1}_{p_j}(\zeta_{p_j}),T^{-1}_{p_j}(z_{p_j})\big)\Big|\Big|\phi\big(T_{p_j}^{-1}(\zeta_{p_j})\big)\Big|\,dS(\zeta_{p_j},z_{p_j})}}\\
+&\underset{II}{\underbrace{\iint_{(z_{p_j},\zeta_{p_j})\in \bd\Om^\epsilon_{p_j}\times(B(0,\delta)\cap  \bd\Om_{p_j})}\Big |H\big(T^{-1}_{p_j}(z_{p_j}),T^{-1}_{p_j}(\zeta_{p_j})\big)\Big|\Big|\phi\big(T_{p_j}^{-1}(\zeta_{p_j})\big)\Big|\,dS(\zeta_{p_j},z_{p_j})}}\\
\lesssim &\Big\|\phi\left(T_{p_j}^{-1}(\cdot)\right)\Big\|_{L^1(\bd\Om_{p_j})} \approx \no{\phi}_{L^1(\bd\Om)},
\end{split}\end{eqnarray}
where $$\Om_{p_j}=\{\rho_{p_j}(z_{p_j}):=\rho\left(T^{-1}_{p_j}(z_{p_j})\right)=P_{p_j}(z_{p_j,1})+r_{p_j}(z_{p_j})<0\}$$ and $P_{p_j}(z_{p_j,1})=F_{p_j}(|z_{p_j,1}|^2)$ or $P(z_{p_j,1})=F(x_{p_j,1}^2)$ as in 
Section \ref{sec:intro}. 
Although the integrals in \eqref{new2} do not cover the full boundaries, the estimate on the complement is trivial because $H$ and is derivatives are uniformly bounded.
Next, note that 
$$\Phi\left(T^{-1}_{p_j}(\zeta_{p_j}),T^{-1}_{p_j}(z_{p_j})\right)=\Phi_{p_j}(\zeta_{p_j},z_{p_j}),$$
where $\Phi_{p_j}$ is the support function of $\Om_{p_j}$, and we can therefore estimate 
\begin{equation}
\begin{split}
\Big |H\big(T^{-1}_{p_j}(\zeta_{p_j}),T^{-1}_{p_j}(z_{p_j})\big)\Big|
\lesssim\frac{1}{|\Phi_{p_j}(\zeta_{p_j},z_{p_j})||\zeta_{p_j}-z_{p_j}|}
\end{split}
\end{equation}
for any $\zeta_{p_j}, z_{p_j}\in \bd\Om$. Nearby each point $p_j$, we will consider domains $\Om_{p_j}$ defined by either \eqref{OmC} or  \eqref{OmR}. 
Here and in what follows, we abuse notation slightly and omit the subscript $p_j$ as well as writing $\phi(\cdot)$ for $\phi\big(T_{p_j}^{-1}(\cdot)\big)$. \\

\noindent{\bf Setting 1: $\Om$ is defined by \eqref{OmC}.} 
We start our estimate of $(I)$ from \eqref{new2} by decomposing the domain of integration and applying  Lemma~\ref{lm22} with $k=1$ to obtain
\begin{equation}
\begin{split}
I=&\iint_{(\zeta,z)\in (\bd\Om\cap B(0,\delta))\times \bd\Om_\epsilon}\cdots\\
=&\iint_{(\zeta,z)\in (\bd\Om\cap B(0,\delta))\times \bd\Om\epsilon \T{~and~} |\zeta_1|\ge|z_1-\zeta_1|}\cdots\\
&+\iint_{(\zeta,z)\in (\bd\Om\cap B(0,\delta))\times \bd\Om\epsilon,|\zeta_1|\le|z_1-\zeta_1| \T{~and~} |z_1|\le 2\delta}\cdots+\iint_{(\zeta,z)\in (\bd\Om\cap B(0,\delta))\times \bd\Om_\epsilon,|\zeta_1|\le|z_1-\zeta_1|  \T{~and~}|z_1| \ge 2\delta}\cdots\\\\
\lesssim&\ (A)+(B)+(C),
\end{split}
\end{equation}
where
\begin{equation}\begin{split}
(A):=&\iint_{(\zeta,z)\in (B(0,2\delta)\cap \bd\Om)\times (B(0,2\delta)\cap \bd\Om_\epsilon)  }\dfrac{|\phi(\zeta)|dS(\zeta,z)}{(\epsilon+|\Im\Phi(\zeta,z)|+F(|z_1-\zeta_1|^2))|z_1-\zeta_1|};\\
(B):=&\iint_{(\zeta,z)\in (B(0,2\delta)\cap \bd\Om)\times (B(0,2\delta)\cap \bd\Om_\epsilon) }\dfrac{|\phi(\zeta)|dS(\zeta,z)}{(\epsilon+|\Im\Phi(\zeta,z)|+F(\frac{1}{2}|z_1|^2))|z_1| };\\
(C):=&\iint_{(\zeta,z)\in (B(0,\delta)\cap \bd\Om)\times \bd\Om_\epsilon \T{~and~}|z_1|\ge 2\delta}\dfrac{|\phi(\zeta)|dS(\zeta,z)}{(\epsilon+|\Im\Phi(\zeta,z)|+F(\frac{1}{2}|z_1|^2))|z_1|}.
\end{split}\end{equation}
$F$ is increasing, so it easily follows that $(C)\lesssim (F(2\delta^2)\delta)^{-1}\no{\phi}_{L^1(\bd\Om)}$.
For (A), we make the change of variables $(\alpha,w)=(\alpha_1,\alpha_2, w_1,w_2)=(\zeta_1,\zeta_2, z_1-\zeta_1, \rho(z)+i\Im\Phi(\zeta,z))$. 
 A direct calculation then establishes that if $\delta$ is chosen sufficiently small  then the Jacobian of this transform does not vanish on 
the domain of integration. Since $\Phi$ is smooth, we can assume that there exists $\delta'>0$
that depends on $\Omega$, $\delta$, and $\rho$ so that if integrate $w_1$ in polar coordinates,
\begin{align*}
(A)&\lesssim \no{\phi}_{L^1(\bd\Om)}\int_{0}^{\delta'}\int_{0}^{\delta'}\frac{r}{(|\Im w_2|+F(r^2))r}\ dr \, d\Im w_2\\
&\lesssim  \no{\phi}_{L^1(\bd\Om)}\int_0^{\delta'}\log F(r^2)\, dr<\infty.
\end{align*}
That the integral is finite follows by the hypotheses on $\phi$ and $F$.\\

Repeating this argument with the change of variables $(\alpha,w)=(\alpha_1,\alpha_2, w_1,w_2)=(\zeta_1,\zeta_2, \frac{1}{\sqrt{2}}z_1, \rho(z)+i\Im\Phi(\zeta,z))$ for the integral (B), we can obtain the same conclusion.\\

To estimate $(II)$ in \eqref{new2}, we use Lemma~\ref{lm22} with $k=1$ and to show the interchanging of $\zeta$ and $z$ is benign. 
In then follows by the same argument as for $(I)$, we obtain $(II)\le \no{\phi}_{L^1(\bd\Om)}$. Therefore, the estimate in 
Seting 1 is complete.\\

\noindent{\bf Setting 2: $\Om$ is defined by \eqref{OmR}.} We omit the proof because it is analogous to Setting 1 with Lemma~\ref{lm23}
replacing Lemma~\ref{lm22}. For  details, see Section 3.2 in \cite{HaKhRa14}.\\

\noindent{\it Part II: Proof of $\no{u}_{L^\infty(\bd\Om)}\lesssim \no{\phi}_{L^\infty(\bd\Om)}$.} The proof of this part is similar to, but simpler than, the argument for  Part III, so we omit it. \\

\noindent{\it Part III: Proof of $\no{u}_{\Lambda^f(\bd\Om)}\lesssim \no{\phi}_{L^\infty(\bd\Om)}$.} We need a general Hardy-Littwood type lemma to prove $f$-H\"older estimates on the boundary.
\begin{lemma} \label{HL}Let $\Om$ be a bounded Lipschitz domain in $\R^N$ and let 
 $G:\R^+\to\R^+$ be an increasing function such that $\dfrac{G(t)}{t}$ is decreasing and $\displaystyle\int_0^s\frac{G(t)}{t}dt<\infty$ for $s>0$ small enough.  If $v\in C^1(\C^n\setminus \bd\Om)$ such that
\begin{eqnarray}\label{5.1}
\quad |\nabla v(x\pm s \nu (x))|\lesssim \frac{G(s)}{s} \quad\T{for any } x\in \bd\Om.
\end{eqnarray}
Then $v\in \Lambda^f(\bd\Om)$ where $f(s^{-1})=\Big(\displaystyle\int_0^s\frac{G(t)}{t}dt\Big)^{-1}.$
\end{lemma}
The proof is basically identical to the corresponding result for domains. See \cite[Theorem 5.1]{Kha13}) for details. Consequently, the focus is now to control the gradient of $H^+$ and $H^-$.
\begin{lemma}\label{lm27}
	For $z\in \bd\Om$, we have
	\begin{enumerate}
		\item $\sup_{z\in \bd\Om_s}|\nabla H^+\phi(z)|\lesssim \frac{G(s)}{s}\no{\phi}_{L^\infty(\bd\Om)}$, and 
			\item $\sup_{z\in \bd\Om^s}|\nabla H^-\phi(z)|\lesssim \frac{G(s)}{s}\no{\phi}_{L^\infty(\bd\Om)}$ 
	\end{enumerate}
where $G(s)=\sup_{p\in \bd\Om}\{\sqrt{F_p^*(s)}\}$ if $\Om$ is defined by \eqref{OmC}and $G(s)=\sup_{p\in \bd\Om}\{\sqrt{F_p^*(s)}|\log\sqrt{F_p^*(s)}|\}$ if $\Om$ is defined by \eqref{OmR}. 
\end{lemma}
\begin{proof}
 Khanh has already proved (1) \cite{Kha13}. For the proof of (2), direct calculations show 
 \begin{align*}
 |\nabla H^-\phi(z)|\lesssim& \no{\phi}_{L^\infty(\bd\Om)}\int_{\bd\Om}\left(\frac{1}{|\Phi(z,\zeta)||z-\zeta|^2}+\frac{1}{|\Phi(z,\zeta)|^2|z-\zeta|}\right)dS(\zeta)\\
 \lesssim &\no{\phi}_{L^\infty(\bd\Om)}\int_{\bd\Om}\frac{dS(\zeta)}{|\Phi(z,\zeta)|^2|z-\zeta|}\\
 \end{align*}
  for $z\in \C^2\setminus\bar\Om$ near $\bar\Om$. 
  We choose a covering $\{B(p_j,\delta)\}^N_{j=1}$ of $\bd\Om$ and change coordinates to set $p_j$ to $0$ as in the proof of Theorem~\ref{t2}; thus our proof reduces to showing
 $$L(z):=\int_{\zeta\in \bd\Om\cap B(0,\delta)}\frac{dS}{|\Phi(z,\zeta)|^2|z-\zeta|}\lesssim \frac{G(\rho(z))}{\rho(z)}, \quad z\in \C^n\setminus\bar{\Om}.$$

For Setting 1, we use Lemma~\ref{lm22} with $k=2$ to interchange the roles of $\zeta$ and $z$. We then estimate
\[
L(z)\lesssim \int_{\zeta\in \bd\Om\cap B(0,\delta)\, \T{and}\, |z_1|\ge |z_1-\zeta_1|}\cdots
+\int_{\zeta\in \bd\Om\cap B(0,\delta)\,\T{and}\,|z_1|\le |z_1-\zeta_1|}\cdots\lesssim (D)+(E),
\]
where 
\begin{eqnarray}
\begin{split}
(D)= &\int_{\zeta\in \bd\Om\cap B(0,\delta)}\frac{dS}{\left(\rho(z)^2+|\Im\Phi(z,\zeta)|^2+F^2(|z_1-\zeta_1|^2)\right)|z_1-\zeta_1|},\\
(E)= &\int_{\zeta\in \bd\Om\cap B(0,\delta)}\frac{dS}{\left(\rho(z)^2+|\Im\Phi(z,\zeta)|^2+F^2(\frac{1}{2}|\zeta_1|^2)\right)|\zeta_1|+\left(\rho(z)^2+|\Im\Phi(z,\zeta)|^2+F^2(\frac{1}{2}|z_1|^2)\right)|z_1|}.
\end{split}
\end{eqnarray}

For the integral $(D)$, if $|z_1-\zeta_1|\ge \delta$ then $(D)\lesssim (F^{2}(\delta^2)\delta)^{-1}$; otherwise we make the change of variable $(w,t)=(z_1-\zeta_1, \Im\Phi(z,\zeta))$. 
We can check that the Jacobian of this tranformation is nonzero on the domain of integration $\delta$ is chosen sufficiently small. Thus,
$$(D)\lesssim (F^{2}(\delta^2)\delta)^{-1}+\int_{|w|\le 2\delta}\frac{dw}{(|\rho(z)|+F(|w|^2))|w|}\le C_\delta \frac{\sqrt{F^*(|\rho(z)|)}}{|\rho(z)|}$$
where the last inequality follows by Lemma 3.2 in \cite{Kha13}.\\

For the integral $(E)$, if $|z_1|\ge \delta$ then $(E)\lesssim (F^{2}(\frac{1}{2}\delta^2)\delta)^{-1}$; otherwise 
$$(E)\lesssim\int_{\zeta\in \bd\Om\cap B(0,\delta)}\frac{dS}{\left(\rho(z)^2+|\Im\Phi(z,\zeta)|^2+F^2(\frac{1}{2}|\zeta_1|^2)\right)|\zeta_1|}.$$
In this case, we make the change variable $(w,t)=(\zeta_1,\Im\Phi(z,\zeta))$. The Jacobian of this transformation is also different zero on the domain of integration if $\delta$ is small. 
We thus obtain the desired estimate for $(E)$.\\

The proof for the real case follows by the same argument using Lemma~\ref{lm23} and Lemma~4.1 in \cite{Kha13}. This is complete the proof of Lemma~\ref{lm27}. 
\end{proof}

Lemma \ref{lm27} allows us to apply Lemma~\ref{HL} to $H^+\phi$ and $H^-\phi$ and establish that $H^+\phi,H^-\phi\in \Lambda^f(\bd\Om)$. 
We may now conclude that $u\in \Lambda^f(\bd\Om)$. 
\end{proof}

\subsection{Proof of Theorem~\ref{t2}}
\begin{proof}[Proof of Theorem~\ref{t2}]
Let $\phi=\sum_{j=1}^{2}\phi_j\,d\bar z_j$ be a bounded, $\mathcal C^1$, $\dib$-closed $(0,1)$-form on $\bar\Om$. 
The solution $u$ of the $\dbar$-equation, $\dbar u=\phi$, provided by the Henkin kernel  is given by
\begin{equation}\label{Henkin}
u=T\phi(z)=H\phi(z)+K\phi(z).\end{equation}
where $H\phi=\int_{\zeta\in \bd\Om}H(\zeta,z)\phi(\zeta)\we\om(\zeta)$ and 
\begin{equation}\begin{split}
K\phi(z)=&\frac{1}{4\pi^2}\int_\Om \frac{\phi_1(\zeta)(\bar\zeta_1-\bar z_1)-\phi_2(\zeta)(\bar\zeta_2-\bar z_2)}{|\zeta-z|^4}\om(\bar\zeta)\wedge\om(\zeta)
\end{split}
\end{equation}

 As mentioned in Section~1, the smoothness of $u$ is a consequence of Theorem 3 in \cite{Ran92}. In particular, Range proved 
	$$\no{T\phi}_{\Lambda_s(\Om)}\lesssim \no{\phi}_{\Lambda_s(\Om)} \quad \T{for all $\phi$ with $\dib\phi=0$ and all $s>0$}$$
	holds on any bounded convex domain $\Om$ in $\C^2$ with smooth boundary. Here $\Lambda_s(\Om)$ is the H\"older space of order $s$.  Thus the proof of Theorem~\ref{t2} will be complete if we prove 
	\begin{eqnarray}\label{L1Dbar}
	\no{T\phi}_{L^1(\bd\Om)}\lesssim \no{\phi}_{L^1(\bd\Om)} 
	\end{eqnarray}
	on our setting of $\Om$.
	For $z\in \bd\Om$ and $\dib \phi=0$ in $\Om$, Shaw \cite[pages 412-414]{Sha89} showed that 
	$$K\phi(z)=-\int_{\zeta\in \bd\Om}H(z,\zeta)\phi(\zeta)\we \om(\zeta).$$
	Although Shaw uses the signed distance to the boundary defining function, her argument is essentially formal and holds for any $C^1$ defining function. Thus we have
	$$\forall z\in \bd\Om, \quad u(z)=\int_{\zeta\in \bd\Om}\left(H(\zeta,z)-H(z,\zeta)\right)\phi(\zeta)\we\om(\zeta).$$
	By the same argument to Part I in Section 2, \eqref{L1Dbar} is obtained. 

\end{proof}

\section{Proof of Theorem~\ref{main1}}\label{sec:proof of thm1}
The next two lemmas are modified versions of Lemmas 4.3 and 4.8 in \cite{Sha91h}.
\begin{lemma}\label{lem:Shaw, d-bar lem}
Suppose that $\Om$ is convex and contains the origin. Let $\alpha$ be a positive, $d$-closed, smooth $(1,1)$-form on $\bar\Om$ supported on $\bar\Om \setminus \overline{B(0,r)}$ for some $r>0$. This means
\[
\alpha = \sum_{j,k=1}^2 \alpha_{j\bar k}\, dz_j \wedge d\z_k
\]
where $\alpha_{j\bar k} \in C^\infty(\bar\Om)$ and $\alpha_{j\bar k}\equiv 0$ on $B(0,r)$. Then there exists a $(0,1)$-form $f$ on $\bar\Om$ so that
\begin{enumerate}
\item $\dbar f=0$;
\item $\p f - \dbar \bar f =\alpha$;
\item There exists $c = c(\Om,r)$ such that
\begin{equation}\label{eqn:f bdd by alpha}
\|f\|_{L^1(\bd\Om)} + \|f\|_{L^1(\Om)} \leq c \|\alpha\|_{L^1(\Om)}.
\end{equation}
\end{enumerate}
\end{lemma}

\begin{proof}
Following Rudin \cite[Theorem 17.2.7]{Rud80}, we let $f(z) = \sum_{k=1}^n f_k(z)\, d\z_k$ where
\[
f_k(z) = \sum_{j=1}^n z_j \int_0^1 t\alpha_{j\bar k}(tz)\, dt.
\]
With this choice of $f$, it follows that both  $\p f - \dbar \bar f = \alpha$ and $\dbar f=0$. 

Since $\|\alpha\|_{L^1(\Om)} = \sum_{j,k=1}^2 \|\alpha_{j\bar k}\|_{L^1(\Om)}$, it follows easily that with induced surface area measure
$d\sigma$,
\begin{align*}
\int_{\bd\Om}|f(z)|\, d\sigma(z) 
&\leq \sum_{j,k=1}^2 \int_{\bd\Om} |z_j| \int_0^1 |t\alpha_{j\bar k}(zt)|\, dt\, d\sigma(z)
\leq \frac {c}{r^2}\sum_{j,k=1}^2 \int_{\bd\Om}\int_0^1 t^3 |\alpha_{j\bar k}(zt)|\, dt\, d\sigma(z) \\
&\leq \frac{c}{r^2}\sum_{j,k=1}^2 \int_\Om |\alpha_{j\bar k}(z)|\, dV(z) = \frac c{r^2} \|\alpha\|_{L^1(\Om)},
\end{align*}
where $dV$ is Lebesgue measure on $\C^2$ and $c$ may change from line to line (and also depends on $\dist(\bd\Om,0)$). Additionally, a similar argument also shows $\|f\|_{L^1(\Om)} \leq c \|\alpha\|_{L^1(\Om)}$. In particular, 
with the change of variables $s = t\tau$,
\begin{align*}
\int_\Om |f(z)|\, dV(z)
&= \int_0^1 \int_{\bd\Om} |f(\tau z)|\, d\sigma(z)\, \tau^3 d\tau
\leq \sum_{j,k=1}^2 \int_0^1 \int_{\bd\Om} \Big| \int_0^1 t\tau z_j \alpha_{j\bar k}(t\tau z)\, dt \Big| d\sigma(z)\, \tau^3d\tau  \\
&\leq \sum_{j,k=1}^2 \int_0^1 \int_{\bd\Om} \int_0^\tau  \big| sz_j \alpha_{j\bar k}(sz)\big| \, ds\, d\sigma(z)\, \tau^2 d\tau  \\
&\leq \sum_{j,k=1}^2  \int_{\bd\Om} \int_0^1 \big| sz_j \alpha_{j\bar k}(sz)\big| \, ds\, d\sigma(z) \leq  \frac c{r^2} \|\alpha\|_{L^1(\Om)}
\end{align*}

\end{proof}

\begin{remark}
In \cite{Sha91h}, Shaw requires that $\alpha$ is positive, i.e., $(\alpha_{j\bar k})$ is a positive definite matrix. In this case, basic linear algebra shows that $2|\alpha_{j\bar k}| \leq \alpha_{j\bar j} + \alpha_{k \bar k}$. 
She then estimates the integral only on the diagonal of $\alpha$. Implicit in her computation is Lelong's computation that positive $(1,1)$ currents $\alpha$ must satisfy $\alpha = \bar\alpha$ and $\alpha_{j\bar k} = -\alpha_{k \bar j}$.
In contrast, we solve the Poincare-Lelong equation for general data with no assumption of positivity. However,  our application to the Nevanlinna class argument only involves positive data.
\end{remark}

\begin{lemma}\label{lem:Shaw d dbar lem}
Suppose that $\Om$ is convex and contains the origin. Let $\alpha$ be a  
$d$-closed, smooth $(1,1)$-form on $\bar\Om$ supported on $\bar\Om \setminus \overline{B(0,r)}$ for some $r>0$.
Then there exists a real-valued function $u\in C^\infty(\bar\Om)$ so that
\begin{enumerate}
\item $i\p\dbar u = \alpha$;
\item $ \| u\|_{L^1(\bd\Om)} \leq c \|\alpha\|_{L^1(\Om)}$ for some constant $c = c(r,\Om) >0$ that is independent of $\alpha$ and $u$.
\end{enumerate}
\end{lemma}

\begin{proof} We use Lemma \ref{lem:Shaw, d-bar lem} to establish the existence of a $\dbar$-closed $(0,1)$-form $f$ that satisfies $\p f - \dbar\bar f=\alpha$ and \eqref{eqn:f bdd by alpha}. Since $f$ is $\dbar$-closed
and in $L^1(\bd\Om)$, we can use Theorem \ref{t2} to establish a function $v$ so that $i\dbar v=f$ and satisfies \eqref{L1bOm}. Note then that
\[
\alpha = i\p \dbar v - i\dbar \p \bar v = i\p\dbar(v+\bar v).
\]
It now follows that $u = v+\bar v$ is the desired function.
\end{proof}

We are now ready to prove Theorem \ref{main1}.
\begin{proof}[Proof of Theorem \ref{main1}.]The second Cousin problem can be solved on convex domains, so there exists $h\in H(\Om)$ with zero set $\hat X$. Extend $h$ to $\C^2$ by setting $h(z) \equiv 1$ for
$z\in \C^2\setminus\Om$.
Let $\alpha = \alpha_{\hat X}$ be the positive
$(1,1)$-current on $\C^2$ defined by $\alpha= i\p\dbar \log|h|$. Observe that $\alpha \equiv 0$ off of $\bar\Om$. 
Let $\vp_\ep\in C^\infty_c(\R)$ be an approximation of the identity, in particular, let $\vp \in C^\infty_c(\R)$, $\supp \vp \subset (-1/2,1/2)$, $\int_\R \vp\, dx=1$, and
$\vp_\ep(x) = \ep^{-1}\vp(x/\ep)$. Let $\Om_\ep = \{z\in \Om : \rho(z) < -\ep\}$.

Define $v_\ep \in C^\infty(\C^2)$ by
\[
v_\ep(z) = \int_{\C^2} \log|h(w)| \vp_\ep(|z-w|)\, dV(w).
\]
Then $v_\ep(z) \to \log|h(z)|$ for almost all $z\in\Om$. By the Poincar\'e-Lelong formula \cite[Theorem 5.1.13]{NoOc90},
$\alpha=0$ on $\{z : h(z)\neq 0\}$, an open set. Therefore, there exists $p\in \Om$ and $r>0$ for which $\alpha|_{B(p,2r)}\equiv 0$. 

Set $\alpha_\ep = i\p\dbar v_\ep$. Then $\alpha_\ep \in C^\infty_{1,1}(\bar\Om)$. Since $d\p\dbar=0$, $\alpha_\ep$ is $d$-closed. Note that if $\ep>0$ is small enough, $\vp_\ep|_{B(p,r)}\equiv 0$. 
Therefore, by translating $p\mapsto 0$,
we can apply Lemma \ref{lem:Shaw d dbar lem} to $\alpha_\ep$ (which we shall do without any further comment regarding the support of $\alpha$ of $\alpha_\ep$). 
Also, $\alpha$ is positive, so $\alpha_\ep$ is as well (on $\Om_\ep$) \cite[Lemma 3.2.13]{NoOc90},
and we write
\[
\alpha_\ep = \sum_{j,k=1}^2 \alpha_{jk}^\ep \, dz_j \wedge d\z_k.
\]
Recall that convolution of a distribution with a test function behaves as follows: $\la T * \vp, \psi \ra = \la T, \tilde \vp * \psi \ra$ where $\tilde \vp(x) = \vp(-x)$. This means
\[
\|\alpha_\ep\|_{L^1(\Om)} = \sup_{g,\ \|g\|_{L^\infty(\Om)}\leq 1} \int_\Om \alpha_\ep \wedge g\, dV.
\]
Each integral in the wedge product is an integral of the $\la \beta_\ep, \psi \ra$ where $\beta_\ep = \beta * \vp$ where $\beta$ is a (positive on $\Om_\ep$)) Radon measure built from the components of $\alpha$. All of this means
\begin{align*}
\|\beta_\ep\|_{L^1(\Om)} &= \sup_{\psi,\ \|\psi\|_{L^\infty(\Om)}\leq 1} \la \beta * \vp_\ep, \psi \ra
= \sup_{\psi,\ \|\psi\|_{L^\infty(\Om)}\leq 1} \la \beta, \tilde \vp_\ep* \psi \ra \\
&= \sup_{\psi,\ \|\psi\|_{L^\infty(\Om)}\leq 1} \int_\Om \tilde\vp_\ep * \psi\, d\beta
\leq \sup_{\psi,\ \|\psi\|_{L^\infty(\Om)}\leq 1} \int_\Om \int_\Om \vp_\ep(x-y) \psi(y)\, dy\, d\beta(x)\\
&\leq \sup_{\psi,\ \|\psi\|_{L^\infty(\Om)}\leq 1}\beta(\Om) \|\psi\|_{L^\infty(\Om)} = \beta(\Om).
\end{align*}
The upshot of this calculation is that because $h$ exactly has $\hat X$ as its zero divisor, the finite area of $\hat X$ guarantees the existence of 
a constant $A>0$ so that $\|\alpha_\ep\|_{L^1(\Om)} \leq A$ where the constant $A$ is independent of $\ep$.

Next, each $v_\ep$ is $d$-closed on $\Om$, so we may invoke Lemma \ref{lem:Shaw d dbar lem} to establish the existence of a real-valued 
$u_{\ep} \in C^\infty(\bar\Om)$ that satisfies $i\p\dbar u_{\ep} = \alpha_\ep$  on $\Om$ and
\[
 \|u_\ep\|_{L^1(\bd\Om)} \leq c \|\alpha_\ep\|_{L^1(\Om)}  \leq cA.
\]
Set $g_\ep = u_\ep - v_\ep$. Then $g_\ep$ is a smooth function on $\bar\Om$ and pluriharmonic on $\Om$ since $i\p\dbar u_\ep = \alpha_\ep = i\p\dbar v_\ep$. 
Moreover, for small $\ep>0$, Lemma \ref{lem:g_ep normal} proves that $\{g_\ep\}$ is a normal family of pluriharmonic functions on $\Om$ 
and therefore there exists a subsequence $\ep_k\to 0$ and a pluriharmonic function $g$ so that $g_{\ep_k}\to g$ uniformly on
compact subsets of $\Om$.

Since $g$ is pluriharmonic on $\Om$, there exists $H \in H(\Om)$ so that
$g = \Rre H$. By construction (and the uniform convergence on compacta), $i\p\dbar g=0$. Define
\[
U(z) = \log|h(z)| + g(z) = \log|e^{H(z)}h(z)|.
\]
The proof is complete once we show that
\[
\int_{\bd\Om_s} |U(z)|\, d\sigma_{\bd\Om_s} \leq C
\]
for some $C>0$ and all $s>0$
but this follows by the argument leading to \cite[(6)]{Gru75}.
\end{proof} 

\begin{lemma}\label{lem:g_ep normal} For $\ep>0$ small, the set of pluriharmonic functions $\{g_\ep\}$ from the proof of Theorem \ref{main1} comprises a normal family. Specifically,
there exists $C>0$ so that if $U \subset\subset \Om$, then there exists $C = C(U)$ that does not depend on $\ep$ so that
$|g_{\ep,s}(z)|\le C$.
\end{lemma}

\begin{proof}
Plurisubharmonic functions are in $L^1_\loc(\Om)$ so $v_\ep = \log|h|*\vp_\ep$ satisfies the following inequality: for $K\subset\Om$ compact, there exists $C_K>0$ so that for every $\ep>0$
\[
\|v_\ep\|_{L^1(K)} \leq C_K.
\]
Following Gruman \cite{Gru75}, we let $U \subset \Om$ have compact closure in $\Om$. 
Let $\eta \in C^\infty_c(\Om)$ so that $0 \leq \eta \leq 1$ and $\eta \equiv 1$ on a neighborhood of $\bar U$. Then for $z\in U$,
\begin{align*}
v_\ep(z) &=  \frac{2}{(2\pi)^2} \int_\Om \frac{1}{|z-w|^2} \triangle\big(  \eta v_\ep(w)\big)\, dV(w) \\
&= \frac{2}{2\pi} \int_\Om \frac{\eta(w)}{|z-w|^2}\triangle v_\ep(w) + \frac{v_\ep(w)}{|z-w|^2}\triangle\eta(w) + \frac{1}{|z-w|^2}\Big( \nabla_w \eta \cdot \nabla_{\w} v_\ep + \nabla_{\w}\eta\cdot\nabla_w v_\ep\Big) dV(w)\\
&= \frac{2}{2\pi} \int_\Om\frac{\eta(w)}{|z-w|^2}\triangle v_\ep(w) \, dV(w) + \frac{2}{2\pi}\int_\Om v_\ep(w)\bigg(\frac{\triangle \eta(w)}{|z-w|^2} - \nabla_{\w}\cdot\Big[ \frac{\nabla_w\eta}{|z-w|^2}\Big] 
+ \nabla_w\Big[\frac{\nabla_{\w}\eta}{|z-w|^2}\Big]\bigg)\, dV(w)
\end{align*}
The second integral is bounded by $C_{\eta}\|v_\ep\|_{L^1(\supp\eta)}$ since $|w-z|$ is bounded away from $0$ since $\supp\nabla\eta$ is a positive distance away from $\bar U$. For the first integral, if $\ep$ is small enough, then
$K = \{ \xi : \xi \in \supp\vp_\ep(w-\cdot) \text{ for any }w\in \bar U\}$ is a compact set in $\Om$. This means
\begin{align*}
\Big| \int_\Om\frac{\eta(w)}{|z-w|^2}\triangle v_\ep(w) \, dV(w) \Big|
&= \int_\Om \int_\Om \frac{\eta(w)}{|z-w|^2} \vp_\ep(w-\xi)\, d\alpha(\xi)\, dV(w) \\
&\leq \int_K |\log |h(\xi)|| \int_\Om \vp_\ep(w-\xi)\eta(w)\frac{1}{|z-w|^2}\, dV(w)\, d\alpha(\xi) \\
&\leq C \alpha(\Om) < \infty
\end{align*}
since $\hat X$ has finite area.
We therefore obtain the bound
\[
\|v_\ep\|_{L^\infty(U)} \leq C_\eta
\]
where $C_\eta$ does not depend on $\ep>0$ (assuming that $\ep>0$ is sufficiently small).
The functions $g_\ep$ are pluriharmonic, so by the Poisson Integral Formula,
\[
g_\ep(z) = \int_{\bd\Om} P(z,w) g_\ep(w)\, d\sigma(w) = \int_{\bd\Om} P(z,w) \big( u_\ep(w) + v_\ep(w)\big)\, d\sigma(w).
\]
Since $u_\ep \in L^1(\bd\Om)$ and $z\in U$ so that $|z-w|$ is bounded away from $0$,
\[
\Big|\int_{\bd\Om} P(z,w)  u_\ep(w)\, d\sigma(w)\Big| \leq C_U \|u_\ep\|_{L^1(\bd\Om)}.
\]
Also, recall that for each fixed $z\in\Om$, $P(z,w) = -\frac{\p}{\p\nu_w} G(z,w)$ where $G(z,w)$ is the Green's function for $\Om$, and $G(z,w)=0$ for all $y\in \bd\Om$. By Green's formula
\begin{align*}
\int_{\bd\Om} P(z,w) v_\ep(w)\, d\sigma(w) &= - \int_{\bd\Om} \frac{\p G(z,w)}{\p\nu_w} v_\ep(w)\, d\sigma(w) \\
&=  \int_\Om v_\ep(w) \triangle G(z,w)\, dV(w) - \int_\Om G(z,w) \triangle v_\ep(w)\, dV(w). 
\end{align*}
Recall that $G(z,w)$ is integrable on $\Om$ in $z$ and in $w$. Indeed, $G(z,w)$ blows up like the Newtonian potential (i.e., integrably) and is symmetric in its arguments.
Consequently,   by Folland \cite[Theorem 6.18]{Fol99}, for $z\in U$, there exists $C = C(U)>0$ so that
\[
\Big| \int_{\bd\Om} P(z,w) v_\ep(w)\, d\sigma(w) \Big| \leq C\big( \|v_\ep\|_{L^\infty(U)} + \|\alpha_\ep\|_{L^1(\Om)}\big).
\]
\end{proof}
 
\bibliographystyle{alpha}
\bibliography{mybib}
\end{document}